\documentclass[12pt,reqno]{amsart}
\usepackage{amsfonts,amssymb,amscd,amstext,mathrsfs,color,amsthm,amsmath}
\usepackage[dvipsnames]{xcolor}
\usepackage[english]{babel}
\usepackage{mathtools}

\usepackage{hyperref}
\usepackage{float}
\usepackage{textcomp}

\usepackage{graphicx}
\usepackage{subcaption}
\usepackage[export]{adjustbox}

\input xy
\xyoption{all}

\theoremstyle{plain}
\newtheorem{theorem}{Theorem}[section]
\newtheorem{lemma}[theorem]{Lemma}
\newtheorem{corollary}[theorem]{Corollary}
\newtheorem{proposition}[theorem]{Proposition}

\theoremstyle{definition}
\newtheorem{definition}[theorem]{Definition}
\newtheorem{example}{Example}
\newtheorem{remark}[theorem]{Remark}

\newtheorem{innercustomthm}{Theorem}
\newenvironment{customthm}[1]
  {\renewcommand\theinnercustomthm{#1}\innercustomthm}
  {\endinnercustomthm}

\numberwithin{equation}{section}

\usepackage{pgf,tikz}
\usetikzlibrary{positioning,shapes,decorations.markings}
\usetikzlibrary{arrows}
\usetikzlibrary[patterns]
\usetikzlibrary{cd}

\usepackage{titlesec} 
\titleformat{\section}[block]{\large\scshape\centering}{\thesection.}{1em}{} 
\titleformat{\subsection}[runin]
  {\normalfont\large\bfseries}{\thesubsection}{1em}{}
\titleformat{\subsubsection}[runin]
  {\normalfont\normalsize\bfseries}{\thesubsubsection}{1em}{}

\usepackage[babel]{csquotes}
\usepackage[backend=biber,style=alphabetic]{biblatex}
\title{Hyperbolicity and Quasi-hyperbolicity in Polynomial Diffeomorphisms of ${\mathbb C}^2$}
\author{Eric Bedford}
\address{Eric Bedford: Stony Brook University, Stony Brook NY 11794}
\email{ebedford@math.stonybrook.edu}

\author{Lorenzo Guerini}
\address{Lorenzo Guerini: Korteweg de Vries Institute for Mathematics, University of Amsterdam, Science Park 107, Amsterdam 1090GE}
\email{lorenzo.guerini92@gmail.com}

\author{John Smillie}
\address{John Smillie: Mathematics Institute, Zeeman Building, University of Warwick, Coventry CV4 7AL}
\email{J.Smillie@warwick.ac.uk}

\begin{document}

\maketitle

\section*{Introduction}
Viewed as dynamical systems, polynomial diffeomorphisms of ${\mathbb C}^2$ have features
in common with complex polynomial maps in one variable and the well studied family of H\'enon maps of ${\mathbb R}^2$.
Friedland-Milnor \cite{FM} show that the dynamically interesting polynomial diffeomorphisms are  finite compositions of generalized (complex) 
 H\'enon mappings.   
 
We begin by defining sets based with the basic dichotomy of bounded/unbounded orbits; $K^\pm$  are the sets where the forward/backward orbits are bounded.  Thus $K^\pm$ and $K:=K^+\cap K^-$ are the analogues of the filled Julia sets for  polynomial maps of ${\mathbb C}$.  The sets where the forward/backward iterates are not equicontinuous are given by $J^\pm:=\partial K^\pm$.  These sets are analogues of the Julia set of polynomial maps in ${\mathbb C}$.  The chaotic dynamics takes place inside the set $J:=J^+\cap J^-$.

Another analogue of the Julia set in dimension two is the  boundary $J^*:= \partial_S K\subset J$, where $\partial_S$ denotes the Shilov boundary (in the sense of function algebras).   In \cite{BS3} we showed that $J^*$ is equal to the closure of the set  ${\mathcal S}$  of saddle periodic points.  Additional dynamical properties/characterizations for $J^*$ were obtained in \cite{BS1,BS3} and \cite{BLS}.



There are several reasons to be interested in polynomial diffeomorphisms which are hyperbolic.
Our understanding of chaotic dynamical systems is most complete in the hyperbolic case. It is also interesting to know how the locus of hyperbolic maps sits inside the parameter space and the relation between hyperbolicity and structural stability (see \cite{DL} and \cite{BerD}). 
We take {\it hyperbolicity} to mean that the system is hyperbolic on its chain recurrent set. 
For polynomial diffeomorphisms  this is equivalent to the condition that the set $J$ is a hyperbolic set, though the chain recurrent set may be larger than $J$ (see \cite{BS1}).  Further, it was  shown that $J=J^*$ for hyperbolic maps, and the stable manifolds ${\mathcal W}^s$ give a Riemann surface lamination of $J^+$.  
In fact, \cite[Theorem~8.3]{BS8} characterized hyperbolicity on $J$ in terms of the existence of  transverse Riemann surface laminations of $J^+$ and $J^-$ in a neighborhood of $J$.  We will strengthen this characterization in Theorem~\ref{thm:lam}.

Hyperbolicity involves uniform expansion and contraction, as well as transversality between expanding and contracting directions.  In \cite{BS8} we defined a canonical metric $\Vert\ \Vert_q^\#$ on the unstable space $E^u_q\subset T_q{\mathbb C}^2$ for each saddle $q\in{\mathcal S}$.  A map is said to be {\it quasi-expanding} if $Df$ expands this metric uniformly (independently of $q\in {\mathcal S}$) in the sense that there is a $\kappa>1$ with  $||Df_qv||_{f(q)}^\#\ge \kappa ||v||_q^\#$  for any nonzero $v\in E^u_q$.  For quasi-expanding maps, this extends to $x\in J^*$.  In \cite{BS8} it was shown that every hyperbolic map is quasi-expanding, with $\Vert\cdot\Vert^\#$  equivalent to the Euclidean norm.  A map is {\it quasi-contracting} if its inverse is quasi-expanding. 

We let $B(q,r)$ denote the ball of radius $r$ centered at $q$ and let $W^s_{q,r}$ denote the connected component of $W^s_q\cap B(q,r)$ containing $x$.  A geometric characterization of quasi-expansion is the Proper, Bounded Area Condition:  there exists $r>0$ such that for all saddles $q\in{\mathcal S}$,  (i) $W^u_{q,r}$ is proper, i.e., closed in $B(q,r)$, and (ii) the area of $W^u_{q,r}$ is uniformly bounded.  This means that the degree of local folding of the manifolds $\{W^u_q: q\in{\mathcal S}\}$ near a saddle point $q_0$ will be bounded.  On the other hand, if  $\{W^u_q: q\in{\mathcal S}\}$  is part of a lamination, then there is no local folding.  

Quasi-expansion may be viewed as a 2-dimensional analogue of semi-hyperbolicity for polynomial maps of ${\mathbb C}$ (see the Appendix of \cite{BS8}).  An important motivation for us was the work of  Carleson, Jones and Yoccoz \cite{CJY}, who showed that semi-hyperbolicity was equivalent to a number of geometric and potential-theoretic properties.

A map will be said to be {\it quasi-hyperbolic} if it is both quasi-expanding and quasi-contracting.  Although there is no explicit requirement of transversality between stable and unstable directions, the stable and unstable manifolds exhibit a weak sort of transversality, as formulated in Proposition \ref{prop:inter}.  Quasi-hyperbolic maps share many properties with hyperbolic ones, and we use this to apply hyperbolic methods in more general contexts.  In the non-hyperbolic case, the canonical metric $\Vert\ \Vert ^\#$ may not be equivalent to the euclidean metric on ${\mathbb C}^2$.  However there is a useful filtration of $J^*$ by a finite number of sets $J^*_{m^s,m^u}$, each of which carries a metric $\Vert\ \Vert^{\#,m^s,m^u}$.   This metric  $\Vert\ \Vert^{\#,m^s,m^u}_x$ is locally equivalent to the Euclidean metric for $x\in J_{m^s,m^u}$, but it may blow up as  $x$ approaches a stratum with larger values of $(m^s,m^u)$.  For maximal values of $(m^s,m^u)$, $J^*_{m^s,m^u}$ is compact and  a uniformly hyperbolic set.


In  Theorem \ref{thm:1}  we represent unstable manifolds as injective holomorphic mappings $\xi_x^u:{\mathbb C}\to {\mathbb C}^2$, and we write  $W^u_x:=\xi^u_x({\mathbb C})$, with ${\mathcal W}^u:=\{W^u_x:x\in J^*\}$.  By part $(ii)$ of Theorem \ref{thm:1},  $\{W^u_q:q\in{\mathcal S}\}$ extends continuously to ${\mathcal W}^u$.   Let us recall  the unstable set
$${\mathbb W}^u_x:=\{z\in{\mathbb C}^2: \lim_{n\to\infty}{\rm dist}(f^{-n}(x),f^{-n}(z))=0\}$$   
which, dynamically, is the attracting basin of $x$  for the inverse map $f^{-1}$.  We refer to $W^u_x$ as an unstable manifold although we only know that $W^u_x\subset {\mathbb W}^u_x$ (Proposition \ref{prop:manif}); we do not know whether these two sets always coincide.

\begin{customthm}{1}\label{thm:1} 
Suppose that $f$ is quasi-expanding.  With the notation above, we have
\begin{itemize}
\item[$(i)$]  For each $x\in J^*$, there is an injective holomorphic immersion $\xi^u_x:{\mathbb C}\to J^+\subset{\mathbb C}^2$ such that $W^u_x\subset{\mathbb W}^u_x$.
\item[$(ii)$]  There exists $\hat r>0$ such that if $0<r<\hat r$, then $W^u_{x,r}$ is a (closed) subvariety of $B(x,r)$.  If $x\in J^*$, then for all but finitely many values of $r<\hat r$, the dependence of the closures on $y$, i.e.  $J^*\ni y\mapsto \overline{ W^s_{y,r}}$, is continuous at $x$ in the Hausdorff topology.
\end{itemize}
\end{customthm}

If $f$ is  quasi-hyperbolic, then Theorem \ref{thm:1} holds also for $\xi^s_x$ and $W^s_x$.  We write ${\mathcal W}^{s/u}$ for the families of the stable/unstable manifolds.  These stable/unstable manifolds are smooth, so it makes sense to say that they have transverse or tangential intersection.  With this, we may characterize uniform hyperbolicity.

\begin{customthm}{2}\label{thm:2} 
Suppose that $f$ is quasi-hyperbolic.  Then $f$ is uniformly hyperbolic on $J$ if and only if there is no tangency between ${\mathcal W}^{s}$ and ${\mathcal W}^u$.
\end{customthm}

Theorem \ref{thm:2} was proved earlier in \cite{BSr} for  {real H\'enon maps of maximal entropy}.  For such maps we have $J\subset{\mathbb R}^2$, and by \cite{BS8} such maps are quasi hyperbolic.  One purpose of the present paper is to extend the work of \cite{BSr} from the context of real, maximal entropy to the more general setting of quasi hyperbolicity.
{}

Theorem \ref{thm:1} will be a consequence of Theorems \ref{thm:LyuPet}, \ref{thm:locman} and Proposition \ref{prop:manif}. Theorem \ref{thm:2} will be proved in~\S\ref{sec:char}.
%

\section{Invariant families of parametrized curves}\label{sec:inv}
\cite{BS8} gives several distinct but equivalent ways of defining quasi-expansion.    One of these is the Proper, Bounded Area Condition, which makes no explicit reference to expansion.  There are also a number of other definitions which use  the pluri-complex Green function $G^+$, which is characterized by the properties:  

\begin{itemize}
\item[$(i)$] $G^+$ is continuous on ${\mathbb C}^2$,  $G^+=0$ on $K^+$, and $G^+>0$ on ${\mathbb C}^2-K^+$, 

\item[$(ii)$] $G^+$ is pluri-subharmonic on ${\mathbb C}^2$, and pluri-harmonic on ${\mathbb C}^2-K^+$, 

\item[$(iii)$] $G^+(z)\le \log(||z||+1) +O(1)$, and $\limsup_{z\to\infty} G^+(z)/\log(||z||+1)=1$.
\end{itemize}

\noindent  A convenient formula (see \cite{H}, \cite{FS}, \cite{BS1}) is that $G^+$ is the super-exponential rate of escape of orbits to infinity:  $G^+(z) = \lim_{n\to\infty}({\rm deg}(f))^{-n}\log(||f^n(z)||+1)$ .

One of the equivalent definitions of quasi-expansion concerns the existence of a large normal family of entire curves imbedded in $J^-$.  We will use this as our definition and derive a number of its properties.

\begin{definition}
Let $f$ be a H\'enon map and $X\subset J$ be an invariant set, so that $J^*\subset \overline{X}$. Suppose that through every $x\in X$ passes an \emph{unstable manifold} $W^u_x$, which is the image of an injective holomorphic immersion $\psi_x:\mathbb C\rightarrow J^-$. Assume further that $\psi_x$ satisfies the normalization conditions
\begin{equation}
\label{eq:norm}
\psi_x(0)=x,\qquad \max_{|\zeta|\le 1} G^+\circ \psi_x(\zeta)=1.
\end{equation}
We say that $f$ is \emph{quasi-expanding} on $X$\footnote{Our definition of quasi expansion is slightly different from \cite{BS8} because of the introduction of the set $X$, which allows us to deal more flexibly with the possibility that $J^*\ne J$.  } if the three following conditions are satisfied:
\begin{itemize}
\item[]\textbf{Invariance:} $f(W^u_x)=W^u_{f(x)}$, for every $x\in X$;
\item[]\textbf{Disjointness:} $W^u_x$ and $W^u_y$ are either equal or disjoint, for every $x,y\in X$;
\item[]\textbf{Normality:} The family $\Psi=\{\psi_x\,:\,x\in X\}$ is normal.
\end{itemize}
We say that $f$ is \emph{quasi-contracting} on $X$ if $f^{-1}$ is quasi-expanding on $X$. A map $f$ is \emph{quasi-hyperbolic} if it is both quasi-expanding and quasi-contracting.
\end{definition} 

We note that if $\psi_x$ satisfies $\eqref{eq:norm}$, then so does the ``rotated''  parametrization $\psi_x(e^{ia}\zeta)$ for any $a\in {\mathbb R}$. Normality of $\Psi$ does not depend on the choice of rotation.

We give the two principal examples of invariant sets $X$, which exist for all H\'enon maps. In Example \ref{ex:saddle} below, the curves are the unstable manifolds of saddle points. For Example \ref{ex:rec} all the curves are the same unstable manifold $W^u_q$ but the parameterizations are different.

\begin{example}[Saddle Points]
\label{ex:saddle} 
We let $X:={\mathcal S}$ be the set of periodic points of saddle type.  Suppose that $f^N(p)=p$, and the multipliers of $Df^N$ at $p$ are $\nu_s,\nu_u\in{\mathbb C}$ with $|\nu_s|<1<|\nu_u|$.  For each $p\in {\mathcal S}$, there is an uniformization $\xi_p:{\mathbb C}\to W^u_p$ such that $\xi_p(0)=p$, and $f^N(\xi_p(\zeta))=\xi_p(\nu_u\zeta)$.  Here $W^u_p$ denotes the standard unstable manifold through the saddle point $p$. We may change the parametrization so that $\psi_p(\zeta):= \xi_p(\alpha\zeta)$ satisfies $\eqref{eq:norm}$.  It follows that $\Psi_{\mathcal S}:=\{\psi_p:p\in{\mathcal S}\}$ satisfies the invariance and disjointness conditions in the definition of quasi-expansion.
\end{example}

\begin{example}[Recentered Unstable Manifold]
\label{ex:rec} 
Let $q$ denote any saddle fixed point, and let $\xi_q:{\mathbb C}\to  W^u_q$  be an uniformization of the unstable manifold through $q$.  By \cite{BLS}, $W^u_q\cap W^s_q$ is a dense subset of $J^*$, and we choose $X\subset W^u_q\cap J$ such that $J^*\subset \overline X$.   For  $y\in X$, let $\zeta_y\in{\mathbb C}$ be such that $\xi_q(\zeta_y)=y$.  We may ``re-center'' the parametrization of this curve to the point $y$, i.e., we choose $\alpha\in{\mathbb C}$ so that $\psi_y(\zeta):=\xi_q(\alpha \zeta + \zeta_y)$ satisfies $\eqref{eq:norm}$.  Thus $\Psi_{q\mathcal R}:=\{\psi_y: y\in X\}$ satisfies the invariance and disjointness conditions in the definition of quasi-expansion.
\end{example}

If $f$ is quasi-expanding on $X$ we let $\widehat\Psi$ denote the set of normal limits of sequences in $\Psi$, and for $x\in\overline X$ we set $\widehat \Psi_x:=\{\psi\in\widehat\Psi:\psi(0)=x\}$.  There are quasi-expanding maps $f$ for which  $\psi\in\widehat\Psi$ may fail to be 1-to-1, and there may be $\zeta_0$ where $\psi'(\zeta_0)=0$.  Thus the elements of $\widehat\Psi_x$ may not be essentially unique modulo rotation of parameters.

Given $\psi\in\widehat\Psi$, we define 
$$
\textrm{Ord}(\psi)=\min\{k\in\mathbb N_{>0}\,|\, \psi^{(k)}(0)\neq 0\}.
$$ 
The following definition will be relevant in section \ref{sec:char}.
\begin{definition}
Let $f$ be quasi-expanding on $X$ and consider a sequence of points $x_k\in X$ converging to $x\in X$. We say that the sequence $x_k$ converges to $x$ \textit{with $m$-folding along the unstable direction} if $\psi_{x_k}\to\psi\in\widehat\Psi$ and $\textrm{ord}(\psi)=m$.

If $f$ is quasi-contracting on $X$  we define an analogous notions of convergence with $m$-folding along the stable direction.

\end{definition}
\begin{proposition}
\label{prop:bigfam}
If $f$ is quasi-expanding, then $\widehat\Psi$ satisfies \eqref{eq:norm} and the following disjointness condition
\begin{equation}
\label{eq:disj}
{\rm If\  }\psi_1,\psi_2\in\widehat \Psi, {\rm \ then\  either\ } \psi_1({\mathbb C})=\psi_2({\mathbb C}),{\rm \  or \ } \psi_1({\mathbb C})\cap \psi_2({\mathbb C})=\emptyset
\end{equation}
\end{proposition}

\begin{proof} It is evident that \eqref{eq:norm} must hold.  If \eqref{eq:disj} fails, there are $\zeta_1,\zeta_2\in{\mathbb C}$ and $\psi_1,\psi_2\in\widehat\Psi$ such that $\psi_1(\zeta_1)=\psi_2(\zeta_2)$, but $\psi_1({\mathbb C})\ne\psi_2({\mathbb C})$.  Thus $\tilde x:=\psi_1(\zeta_1)=\psi_2(\zeta_2)$ is an isolated point of $\psi_1({\mathbb C})\cap \psi_2({\mathbb C})$.  On the other hand, $\psi_1$ is the locally uniform limit of $\psi_{1,j}\in\Psi$ (and similarly for $\psi_2$).  By the continuity of complex intersections, there is an intersection point of  $\psi_{1,j}({\mathbb C})\cap \psi_{2}({\mathbb C})$ near $\tilde x$ when $j$ is sufficiently large.  Now if we choose $k$ sufficiently large, there is an intersection point of $\psi_{1,j}({\mathbb C})\cap \psi_{2,k}({\mathbb C})$ near $\tilde x$.  This contradicts disjointness in the definition of quasi-expansion.
\end{proof}

Given $\psi_1,\psi_2\in\widehat\Psi_x$ we have $\psi_1(0)=\psi_2(0)=x$, and therefore $\psi_1(\mathbb C)=\psi_2(\mathbb C)$. This allows us to define an unstable manifold through every point of $\overline X$.
\begin{definition}
Let $f$ be quasi-expanding on $X$. Given $x\in \overline X$ we define $W^u_x:=\psi(\mathbb C)$, where $\psi$ is any element of $\widehat\Psi_x$.
\end{definition}

Since the sets $W^u_x$ are uniform limits of disjoint disks, we have the following consequence of \cite[Proposition 12]{LP}.

\begin{theorem}[Lyubich-Peters]
\label{thm:LyuPet}
The set $W^u_x\subset J^-$ is biholomorphic to $\mathbb C$.
\end{theorem}
\begin{proof}
Every $\psi\in\widehat\Psi_x$ is the locally uniform limit of a sequence of maps whose images lie in $J^-$, therefore $W_x^u=\psi(\mathbb C)\subset J^-$. The map $\psi$ is not constant by \eqref{eq:norm}, therefore by \cite[Proposition 12]{LP}  $W^u_x$ has no singular points. In particular $W^u_x$ is a Riemann surface.

Every $\psi\in \widehat\Psi^s_x$ is finite-to-one, see for instance \cite[Lemma 3.3]{GP}. Therefore $W^u_x$ is simply connected. By Liouville theorem the set $W^u_x\subset \mathbb C^2$ is not biholomorphic to the unit disk nor the Riemann sphere. We conclude that $W^u_x$ is biholomorphic to $\mathbb C$.
\end{proof}

For $x\in\overline X$ and $r>0$ we let $B(x,r)\subset{\mathbb C}^2$ denote the ball of radius $r$, centered at $x$.  For $\psi\in\widehat\Psi_x$, we let ${\mathcal D}_\psi={\mathcal D}_{\psi}(r)$ denote the connected component of the open set $\psi^{-1}(B(x,r))$ containing the origin.  By \eqref{eq:norm}, $\psi$ is non-constant, and we can choose $r>0$ small enough that  $\psi:{\mathcal D}_\psi\to B(x,r)$ is proper, which corresponds to the property that ${\rm dist}(\psi,x)=r>0$ on $\partial{\mathcal D}_\psi$.  We may choose a uniform $r$ that works for all $x\in J^*$.  By the Maximum Principle, ${\mathcal D}_\psi$ is a topological disk, so by the Riemann Mapping Theorem it is conformally equivalent to the unit disk.  For $x\in\overline X$ and $\psi\in\widehat\Psi_x$, the fact that $\psi$ is non constant means that there is an integer $m\ge1$ and $\vec a\in{\mathbb C}^2$ such that $\psi(\zeta)=x + \vec a\zeta^m + \cdots$, where the dots indicate higher powers of $\zeta$.  Thus there are $r,\rho>0$ such that ${\rm dist}(\psi(\zeta),x)\ge r$ for all $|\zeta|=\rho$.  This means that ${\mathcal D}_{\psi}\subset\{|\zeta|<\rho\}$.  Since $|\psi'|$ is bounded by some number $M$ on $|\zeta|\le\rho$, we know that $\{|\zeta|<r/M\}\subset {\mathcal D}_\psi$.  By the compactness of $\widehat\Psi$, we have:

\begin{proposition}
\label{prop:radii}
For sufficiently small $r>0$, there are $0<\rho_1<\rho_2$  such that 
\begin{equation}
\label{eq:disks}
\{|\zeta|<\rho_1\}\subset{\mathcal D}_\psi\subset\{|\zeta|<\rho_2\}\qquad\forall \psi\in\widehat\Psi.
\end{equation}
\end{proposition} 
By {\it inner radius} $\rho^{\rm in}:=\rho^{\rm in}(r)$ and {\it outer radius} $\rho^{\rm out}:=\rho^{\rm out}(r)$ we will denote the maximal/minimal values such that \eqref{eq:disks} holds.  From Proposition \ref{prop:radii}, is evident that both of these radii are monotone increasing in $r$ and tend to zero as $r\to0$.

The converse of the previous Proposition is also valid
\begin{proposition}
\label{prop:radii2}
Given $\rho>0$, there exist $0<r_1<r_2$ such that
$$
\mathcal D_{\psi}(r_1)\subset \{|\zeta|<\rho\}\subset \mathcal D_\psi(r_2).\qquad\forall\psi\in\widehat\Psi.
$$
\end{proposition}
\begin{proof}
By quasi-expansion there exists $M>0$ so that $\sup_{\zeta\in \{|\zeta|<\rho\}}\Vert\psi'(\zeta)\Vert<M$ uniformly on $\widehat\Psi$. The second inclusion follows by taking $r_2=M\rho$. It is not hard to show that for $\psi\in \widehat\Psi$ there exists $r>0$ so that $\mathcal D_\psi(r)\subset \{|\zeta|<\rho\}$. The existence of $r_1$ follows from compactness of $\widehat\Psi$.
\end{proof}

\section{Expanded metric}
\label{sec:exp}
Let $f$ be quasi-hyperbolic on $X$.  If $x\in X$, then $f(\psi_x)$ differs from $\psi_{f(x)}$ by a linear reparametrization, so there exists $\lambda_x\in{\mathbb C}$ such that 
\begin{equation}
\label{eq:conj}
f(\psi_x(\zeta))=\psi_{f(x)}(\lambda_x\zeta).
\end{equation}


The second condition in \eqref{eq:norm} may be interpreted as saying that the disk $\{|\zeta|<1\}$ is the largest disk centered at the origin and contained in the set $\{G^+\circ\psi_x<1\}\subset{\mathbb C}$.  Under $f$ (equivalently $L_{\lambda_x}(\zeta)=\lambda_x\,\zeta$) this unit disk is taken to the disk $\{|\zeta|<|\lambda_x|\}\subset{\mathbb C}$.  A basic property of  $G^+$ is that $G^+\circ f=deg(f)G^+$, and $deg(f)\ge2$.   Since $\{|\zeta|<|\lambda_x|\}$ is then the largest disk inside $\{G^+<{\rm deg}(f)\}$, it follows that $|\lambda_x|>1$.

We define three rate of growth functions: 
$$m_\psi(r):=\max_{|\zeta|\le r} G^+(\psi(\zeta)), \quad m(r): =\inf_{\psi\in\widehat\Psi}m_\psi(r), \quad M(r): = \sup_{\psi\in\widehat\Psi}m_\psi(r)$$
With this notation, we see that for $x\in X$, the value of $|\lambda_x|$ is defined by the condition
$ m_{\psi_{f(x)}}(|\lambda_x|) = d$.

\begin{proposition}
\label{prop:eqcond}
If $f$ is quasi-expanding, then the following hold:
\begin{itemize}
\item[(i)]  $M(r)<\infty$ for all $r<\infty$.
\item[(ii)]  $m(r)>0$ for all $r>0$.
\item[(iii)]  Let $\kappa$ be such that  $M(\kappa)={\rm deg}(f)$.  Then $\kappa>1$, and $|\lambda_x|\ge \kappa$ for all $x\in X$.
\item[(iv)] There exists $K\ge \kappa$ so that $|\lambda_x|\le K$ for all $x\in X$.
\end{itemize}
\end{proposition}

\begin{proof}  
Properties (i-iii) are proved in \cite{BS8}. Property (iv) is shown in \cite[Corollary 2.8]{GP}.  
\end{proof}

In fact, each of the conditions (i-iii) is equivalent to normality of $\Psi$ (see \cite{BS8}).

Given a quasi-expanding map, it is possible to define (intrinsic) metrics and distances on each unstable manifold $W^u_x$. With respect to these quantities, the map $f$ expands uniformly. These definitions are given in term of the injective immersions $\psi_x$, but do not depend on the particular choice of $\psi_x$ (recall that $\psi_x$ is uniquely determined up to rotations).

\subsubsection*{Expanding metric:} For each $x\in X$, we define $E^u_x$ to be the 1-dimensional subspace of $T_x({\mathbb C}^2)$ spanned by $\psi'_x(0)$, or equivalently the tangent space of $W^u_x$ at the point $x$.  We let $|\cdot|_e$ denote the Euclidean norm and define a new norm on $E^u_x$:
$$\Vert v \Vert^\#_x := \frac{| v |_e}{|\psi'_x(0)|_e}, \qquad v\in E^u_x. $$
With this definition, the (operator) norm of $\psi'_x(0)$ is 1. By applying the chain rule on $f\circ\psi_x$, we find that the operator norm of the restriction of $D_xf$ to $E^u_x$ with respect to the metric $\Vert\cdot\Vert^\#_x$ is given by
$$\Vert D_xf\Vert_x^\# = \frac{\Vert D_xf(v)\Vert ^\#_{f(x)}} {\Vert v\Vert ^\#_x} = |\lambda_x|, \qquad {\rm\ for\ }v\in E^u_x$$

\subsubsection*{Expanding distance:}For each $x\in X$, we define $\mathrm{dist}^u_x\,:\,W^u_x\times W^u_x\rightarrow \mathbb R_+$, by pulling back the euclidean distance of $\mathbb C$ via $(\psi^u_x)^{-1}$. We have that
$$
{\rm dist}^u_{f(x)}\left(f(y_1),f(y_2)\right)=|\lambda_x|\,{\rm dist}^u_x(y_1,y_2),\qquad {\rm\ for\ }y_1,y_2\in W^u_x
$$
If $x,y\in X$ are so that $W^u_y=W^u_x$, then the holomorphic maps $\psi^u_x$ and $\psi^u_y$ coincide up to an affine transformation. Therefore $\mathrm{dist}^u_y=c \,\mathrm{dist}^u_x$ for some real number $c>0$. 
\begin{lemma}
\label{intrinsicmetric}
Suppose that $W^u_x=W^u_y$ for some $x,y\in X$. Given $a\in \mathbb C^*$ and $b\in\mathbb C$ so that $\psi_y(\zeta)=\psi_x(a\zeta+b)$, we have that
$$
||a|-1|\le |b|=\mathrm{dist}^u_x(x,y).
$$
If $c>0$ is so that $\mathrm{dist}^u_y=c\,\mathrm{dist}^u_x$, then $|c^{-1}-1|\le \mathrm{dist}^u_x(x,y)$.
\end{lemma}
\begin{proof}
The first normalization condition in \eqref{eq:norm} implies $y=\psi_y(0)=\psi_x(b)$, showing that $|b|=\mathrm{dist}^u_x(x,y)$. By the second condition in \eqref{eq:norm}, there exists $|\zeta_y|=1$, so that 
$$
1=G^+\circ\psi_y(\zeta_y)=G^+\circ\psi_x(a\zeta_y+b),
$$
proving that $|a\zeta_y+b|\ge 1$, and therefore that $|a|\ge 1-|b|$. Similarly we may find $|\zeta_x|=1$ so that 
$$
1=G^+\circ\psi_x(\zeta_x)=G^+\circ\psi_y\left(\frac{\zeta_x-b}{a}\right),
$$
proving that $\left|\frac{\zeta_x-b}{a}\right|\ge 1$, and therefore $|a|\le 1+|b|$. It is not difficult to show that $c=|a|^{-1}$, concluding the proof of the lemma.
\end{proof}

Properties (iii-iv) of Proposition \ref{prop:eqcond} give the following
\begin{proposition}
\label{expansion2}
Let $f$ be quasi-expanding. Then there exists $1<\kappa\le K$ so that for every $x\in X$ we have $$\kappa\le \Vert D_xf\Vert_x^\#\le K,$$and
$$
\kappa\, \mathrm{dist}^u_x(y_1,y_2)\le \mathrm{dist}^u_{f(x)}(f(y_1),f(y_2))\le K\, \mathrm{dist}^u_x(y_1,y_2)\qquad {\rm\ for\ }y_1,y_2\in W^u_x,
$$
\end{proposition}

Let $r>0$ be sufficiently small, and write $\rho^{\mathrm{in}},\rho^{\mathrm{out}}$ for the inner and outer radii given by Proposition \ref{prop:radii}. As in the introduction, $W^u_{x,r}$ denotes the connected component of $W^u_x\cap B(x,r)$ containing a given $x\in\overline X$. Given $\psi\in\widehat\Psi_x$ we have $W^u_{x,r}=\psi(\mathcal D_\psi)$. We denote the family of these disks by ${\mathcal W}^u_r:=\{W^u_{x,r}:x\in\overline X\}$.

\begin{proposition}
\label{prop:contr} 
If $f$ is quasi-expanding, then for sufficiently large $N$, $f^{-N}$ maps ${\mathcal W}^u_r$ inside itself.  That is, for each $x\in\overline X$, $f^{-N}(W^u_{x,r})\subset W^u_{{f^{-N}(x)},r}$.  Further, if $y\in W^u_{x,r}$, then $dist(f^{-n}(x),f^{-n}(y))\to0$ like a constant times $\kappa^{-n}$ as $n\to +\infty$.
\end{proposition}

\begin{proof}   
Let $\kappa$ is as in Proposition \ref{prop:eqcond}, and let $N$ be  large enough that $\rho_1\kappa^N>\rho_2$.  Thus  
$$f^{-N}W^u_{x,r} \subset  \psi_{f^{-N}(x)}(|\zeta|<\kappa^{-N} \rho^{\mathrm{out}}) \subset \psi_{f^{-N}(x)}({\mathcal D}_{f^{-N}(x)})= W^u_{{f^{-N}(x)} ,r}$$
For the last assertion, by the normality of $\Psi$ we have $|\psi'(\zeta)|\le C$ for all $\psi\in\Psi$ and $|\zeta|\le \rho^{\mathrm{out}}$.  Now for any two points of $W^u_x$, we may write them as $y'=\psi_x(\zeta')$ and $y''=\psi_x(\zeta)$.  The distance between the $\zeta$-coordinates of $f^{-n}(y')$ and $f^{-n}(y'')$ is $\kappa^{-n}|\zeta'-\zeta''|$.  Thus we conclude that ${\rm dist}(f^{-n}(y'),f^{-n}(y'')) \le \rho^{\mathrm{out}}\kappa^{-n}C$, which proves the last assertion. 
\end{proof}

\begin{proposition} 
\label{prop:manif} 
If $f$ is quasi-expanding, then for all $x\in \overline X$, $W^u_x =\bigcup_{n\ge 0} f^{n}W^u_{ {f^{-n}(x)},r}$  is a manifold consisting of unstable points, i.e. $W^u_x\subset{\mathbb W}^u_x$.
\end{proposition}

\begin{proof} 
Given $x\in \overline X$ we choose a sequence $\psi_n\in\widehat\Psi^s_{f^n(x)}$. Then $W^u_{{f^{-n}(x)}, r}\supset \psi_n(|\zeta|<\rho^{\mathrm{in}})$.  By quasi-contraction, we have that $f^{n}(W^s_ {{f^{-n}(x)},r} )\supset \psi_0(|\zeta|<\kappa^n\rho^{\mathrm{in}})$. It follows that
$$
W^u_x=\psi_0(\mathbb C)=\bigcup_{n\ge 0} f^{n}W^u_{ {f^{-n}(x)},r}.
$$
Proposition \ref{prop:contr} shows that $W^u_{x,r}$, and thus $W^u_x$, is contained in ${\mathbb W}^u_x$.
\end{proof}

Knowing that $W^u_x\subset \mathbb W^u_x$, we can show that the family of curves $\mathcal W^u=\{W^u_x\,:\,x\in \overline X\}$ satisfies the three condition of quasi-expansion. This can be done as in section 3 of \cite{GP}\footnote{In that paper it is assumed that the set $X\subset J^*$, but the proof also works in our case}. This shows the following
\begin{theorem}
\label{thm:closure}
Let $X\subset J$ be an $f$-invariant set so that $J^*\subset \overline{X}$. Then $f$ is quasi-expanding on $X$ if and only if $f$ is quasi-expanding on $\overline X$.
\end{theorem}
\begin{remark}
If $f$ is quasi-expanding on the closed and invariant set $X$, then by \cite[Theorem 3.5]{GP} the unstable manifolds are uniquely determined. More precisely, given two families of unstable manifolds $\mathcal W^{u,1}$ and $\mathcal W^{u,2}$  for which $f$ is quasi-expanding on $X$, then  $W^{u,1}_x=W^{u,2}_x$ for every $x\in X$.
\end{remark}

We conclude this section with a uniform Lyapunov stability condition within leaves of ${\mathcal W}^u$.

\begin{proposition} 
If $f$ is quasi-expanding, and if $\epsilon>0$ is given, then there exists $\delta>0$ such that if $x,y\in X$, $dist(x,y)<\delta$, and $y\in W^u_{x,r}$, then $dist(f^{-n}(x),f^{-n}(y))<\epsilon$ for all $n\ge0$.
\end{proposition}

\begin{proof}  Let $\kappa$ as in Proposition \ref{prop:eqcond}.  By \cite[Theorem 3.2]{BS8}, the inner/outer radii in \eqref{eq:disks} have a ratio $a:=\rho^{\rm out}/\rho^{\rm in}$ which is determined by the local area bound of ${\mathcal W}^u$.  The sets $W^u_{x,r}$, and thus the area bound, are unchanged if we replace $f$ by $f^n$.  Thus we may assume that $\kappa>a$.  We let $M$ denote the supremum of the euclidean norm $\Vert Df^{-1}(x)\Vert$ for $x$ in an $\epsilon$ neighborhood of $X$.  Now let us choose $\delta$ such that $\rho^{\rm out}(\delta)<\rho^{\rm in}(\epsilon/M^n)$.  It follows that the image $f^{-n}(W^u_{x,\delta})\subset W^u_{f^{-n}(x),\delta}$ since 
$$|\lambda_x\lambda_{f^{-1}(x)}\cdots\lambda_{f^{-n+1}(x)}| {\mathcal D}_{\psi_{x}} (\delta)\subset \kappa^{-n}{\mathcal D}_{\psi_{x}}\subset {\mathcal D}_{\psi_{f^{-n}(x)}}(\delta)$$
By the choice of $M$, we see that $f^{-j}(W^u_{x,\delta})\subset W^u_{f^{-j}(x),\epsilon}$ for $0<j<n$. 
\end{proof}

\section{Motion of local manifolds}
\label{sec:locunst}
Let $f$ be quasi-expanding on an invariant set $X$. By Theorem \ref{thm:closure} we may assume that $X$ is compact. Throughout this section we will assume that $\hat{r}$ is small enough for Proposition \ref{prop:radii}.  Recall that for every $r\le\hat r$ and $\psi\in\widehat\Psi_x$ the map $\psi:\mathcal D_\psi(r)\rightarrow W^u_{x,r}$ is proper. Let $\rho^{\rm in/out}:=\rho^{\rm in/out}(\hat r)$ as in Proposition \ref{prop:radii}.

\begin{proposition} 
There exists $C<\infty$ such that for all $x\in X$ and balls $B(x,r)\subset B(x,\hat r)$, the number of connected components of $W^u_{x,\hat r}\cap B(x,r)$ is less than $C$.
\end{proposition}

\begin{proof}
Let $\Omega:=\{\zeta\in{\mathcal D}_\psi(\hat r):\psi_x(\zeta)\in B(x,r)\}$, so the components of $\Omega$ are the pre-images of the components of $ W^u_{x,r}\cap B(x,\rho)$.  It follows that $\sigma(\zeta):={\rm dist}^2(\psi_x(\zeta),x)$ gives a proper map $\sigma:\Omega\to[0,\rho^2)$.  Since $\psi_x$ is holomorphic, it follows from the maximum principle that all components of $\Omega$ are simply connected.  Thus the number of connected components is bounded by the total number of critical points of $\sigma$ (counted with multiplicity) inside $\Omega\subset{\mathcal D}_\psi(\hat r)\subset \{|\zeta|<\rho^{\rm out}\}$.  The number of critical points is bounded since $\Psi$ is a normal family.
\end{proof}

We define $\widehat W^u_{x,r}:= \bigcap_{\epsilon>0} W^u_{x,r+\epsilon} \cap B(x,r)$.  We note that $\widehat W^u_{x,r}$ will consist of $W^u_{x,r}$, together with possibly a finite number of connected components of $W^u_x\cap B(x,r)$ which become ``joined'' to $W^u_{x,r}$ at critical points inside $B(x,r+\epsilon)$ for every $\epsilon>0$.

Recall that the \emph{Hausdorff distance} between two compact sets $X,Y\subset \mathbb C^2$ is given by $d_H(X,Y)=\max\{\partial_H(X,Y),\partial_H(Y,X)\}$, where $$\partial_H(X,Y)=\sup_{x\in X}\inf_{y\in Y} {\rm dist}(x,y).$$ 
\begin{proposition}
\label{prop:convergencelocal}
Let $x_k\in X$ be a sequence converging to $x\in X$, and $r_k>0$ be a sequence converging to $r\in (0,\hat r)$. Then
$$
\lim_{k\to\infty}\partial_H\left(\overline{W^u_{x,r}},\overline{W^u_{x_k,r_k}}\right)= 0\quad\text{and}\quad\lim_{k\to\infty}\partial_H\left(\overline{W^u_{x_k,r_k}},\overline{\widehat W^u_{x,r}}\right)= 0
$$
\end{proposition}
\begin{proof}
It is sufficient to consider sequences for which $r_k\in (0,\hat r)$, and for which $\psi_{x_k}$ converges locally uniformly to some $\psi\in\widehat\Psi_{x}$. 

{\bf Part 1} ($\lim_{k\to\infty}\partial_H\left(\overline{W^u_{x,r}},\overline{W^u_{x_k,r_k}}\right)= 0$). The map $s\mapsto \overline{W^u_{x,s}}$ is always continuous in a left neighborhood of the point $r$. Given $\varepsilon>0$, we choose $0<\delta<\varepsilon$ so that $\partial_H(\overline{W^u_{x,r}},\overline{W^u_{x,r-3\delta}})<\varepsilon$. For $k$ sufficiently large, we have $r_k\ge r-\delta$, which implies that $\partial_H(\overline{W^u_{x_k,r-\delta}},\overline{W^u_{x_k,r_k}})=0$. Part 1 follows by the triangle inequality, once we show that for $k$ large
\begin{equation}
\label{eq:part1}
\partial_H(\overline{W^u_{x,r-3\delta}},\overline{W^u_{x_k,r-\delta}})<\delta.
\end{equation}

Let $D$ be the connected component of $\psi^{-1}(\overline{W^u_{x,r-3\delta}})$ containing the origin. The set $D$ is closed and contained in $\{|\zeta|<\rho^{\rm out}\}$, thus it is compact. For $k$ large we have $\sup_{\zeta\in D}{\rm dist}(\psi(\zeta),\psi_k(\zeta))<\delta$, and therefore $\partial_H(\overline{W^u_{x,r-3\delta}},\psi_k(D))<\delta$. For $k$ large we have $B(x,r-3\delta)\subset B(x_k,r-2\delta)$, showing that $\psi_{x_k}(D)$ lies within the ball $B(x_k,r-\delta)$. This implies that $\psi_{x_k}(D)\subset\overline{W^u_{x_k,r-\delta}}$, and therefore we obtain \eqref{eq:part1}.

{\bf Part 2} ($\lim_{k\to\infty}\partial_H\left(\overline{W^u_{x_k,r_k}},\overline{\widehat W^u_{x,r}}\right)= 0$). Let $y_k\in \overline{W^u_{x_k,r_k}}$ be a point lying at the furthest distance from the set $\overline{\widehat W^u_{x,r}}$, and let $\zeta_k=\psi_{x_k}^{-1}(y_k)$. Since $\overline{\widehat W^u_{x,r}}\subset W^u_{x,\hat r} $, by Proposition \ref{prop:radii} we have $|\zeta_k|<\rho^{\rm out}$. By taking a subsequence if necessary, we may therefore assume that $\zeta_k$ converges to a point $\zeta_\infty$. By local uniform convergence, the sequence $y_k=\psi_{x_k}(\zeta_k)$ converges to the point $y=\psi(\zeta_\infty)\in W^u_{x}$. Part 2 follows once we show that $y\in \overline{\widehat W^u_{x,r}}$. 

The norm of the derivative $\Vert\psi'\Vert$ is uniformly bounded from above on the disk $\{|\zeta|<\rho^{\rm out}\}$ for all $\psi\in\widehat\Psi$. Therefore for every $\delta>0$ sufficiently small, and $k$ large so that $\zeta_k$ and $\zeta_\infty$ are close enough, we may assume that $\psi_{x_k}([\zeta_k,\zeta_\infty])\subset B(x_k,r+\delta)$. It follows that $\psi_{x_k}(\zeta_\infty)\in W^u_{x_k,r+\delta}$, or equivalently that $\zeta_\infty\in\mathcal D_{\psi_{x_k}}(r+\delta)$.

Let $\gamma_k:[0,1]\rightarrow \mathcal D_{\psi_{x_k}}(r+\delta)$ be a continuous curve connecting the point $0$ to $\zeta_\infty$. Assume that $k$ is sufficiently large so that $\sup_{|\zeta|<\rho^{\rm out}}{\rm dist}(\psi(\zeta),\psi_{x_k}(\zeta))<\delta$. Then the curve $\psi\circ \gamma_k$ lies within the ball $B(x,r+2\delta)$. This curve connects the points $x=\psi\circ\gamma_k(0)$ and $y=\psi\circ\gamma_k(1)$, and remains inside the unstable manifold $W^u_x$. Therefore $y\in W^u_{x,r+2\delta}$ for every $\delta>0$. This implies that $y\in \overline{\widehat W^u_{x,r}}$, concluding the proof of the Proposition.
\end{proof}

\begin{definition}
We say that $0<r<\hat r$ is a {\it regular radius} for $x\in X$ if $\overline{W^u_{x,r}}$ intersects $\partial B(x,r)$ transversally.
\end{definition}


\begin{proposition}
Given $x\in X$, there are only finitely many values $r$ which are less than $\hat r$ and which are not regular for $x$. 
\end{proposition}
\begin{proof}
Since $W^u_{x}$ is the image of $\psi_{x}$, having a tangential intersection between $\psi_{x}({\mathbb C})$ and  $\partial B(x,r)$ at a point $\psi_{x}(\zeta_x)$, $\zeta_x\in{\mathcal D}_\psi(\hat r)$ means that $\zeta_x$ is a critical point for the map $\zeta\mapsto{\rm dist}^2(\psi_{x}(\zeta),x)$.  Since this distance map is the square of the modulus of a holomorphic mapping, there can be only finitely many critical points in a compact region.  The (finitely many) critical values of this function correspond to the squares of the values of $r$ which are not regular.   
\end{proof}

\begin{theorem} 
\label{thm:locman}
Let $f$ be quasi-expanding on the compact invariant set $X$, and let $r<\widehat r$. The map $\psi:{\mathcal D}_\psi(r)\to B(x,r)$ is a proper imbedding.  Further, $W^u_{x,r}$ has the following  properties:
\begin{itemize}
\item[(i)]  $W^u_{x,r}$ is a nonsingular subvariety of $B(x,r)$.  
\item[(ii)]  The family ${\mathcal W}^u_r$ has uniformly bounded area:  $\sup_{x\in X}{\rm Area}(W^u_{x,r})<\infty$.
\item[(iii)] If $r$ is a regular radius for $x$, then the function $(y,s)\mapsto \overline{W^u_{y,s}}$ is continuous in the Hausdorff topology at $(x,r)$.
\end{itemize}
\end{theorem}

\begin{proof}  

Item (i) is a consequence of Theorem \ref{thm:LyuPet}.   Item (ii)  is an easy consequence of the normality of $\Psi$. If $r$ is regular for $x$ then the map $s\mapsto \overline{W^u_{x,s}}$ is continuous with respect to the Hausdorff distance at $r$. Therefore we have that $\widehat W^u_{x,r}=W^u_{x,r}$. Point (iii) is a consequence of Proposition \ref{prop:convergencelocal}.
\end{proof}


Next we observe that  if $f$ is quasi-hyperbolic, the local disks $W^s_{x,r}$ and $W^u_{x,r}$ have a weak sort of transversality, in the sense that if two centers are close, then the local stable/unstable disks themselves must intersect. 

\begin{proposition}
\label{prop:inter}
Let $f$ be quasi-hyperbolic on $X$, and $r<\hat r$.  Then for $r_1<r$ there exists $r_0>0$ with the following property.  If $B_0\subset B_1\subset{\mathbb C}^2$ are concentric balls such that $B_j$ has radius $r_j$, then for $x',x''\in B_0\cap X$, the intersection   $W_{x',r}^s \cap  W_{x'',r}^u \cap B_1$ is nonempty.
\end{proposition}

\begin{proof} 
By Theorem \ref{thm:locman} we know that if we fix $x$ and change $r$ slightly, $X\ni y\mapsto W^{s/u}_{y,r}$ is continuous at $x$ with respect to the Hausdorff topology.  For $r_0>0$ sufficiently small, we have $B(x',r_1)\subset B(x'',r)$ for all $x',x''\in B_0$.  Thus for $x\in B_0$,  $W^{s/u}_{x,r}\cap B_1$ are closed subvarieties of $B_1$.   Since these varieties do not actually coincide in a neighborhood of $x$,  $W^u_{x,r}\cap W^s_{x,r}=\{x\}$ in a neighborhood of $x$.  The intersection of complex varieties has the following continuity property:  For any $r_1>0$, we may choose $r_0>0$ sufficiently small that for $x',x''\in B_0$,  $W^u_{x',r}\cap W^s_{x'',r}\cap B_1\ne\emptyset$. 
\end{proof}

We conclude this section with a result that is proved much along the lines of Theorem~9.6 of \cite{BLS}:

\begin{theorem}
\label{thm:dens}
Let $f$ be quasi-hyperbolic, and let $T_1,T_2\subset J^*$ be closed, invariant subsets.   Then $\bigcup_{x\in T_1} W^s_x\cap J^*$ and $(\bigcup_{x\in T_1}W^s_x)\cap (\bigcup_{y\in T_2}W^u_y)$ are dense  subsets of $J^*$.
\end{theorem}

\section{Projection maps}
Let $f$ be quasi-expanding on $X$. Given $x\in X$ we define $P_x:\mathbb C\rightarrow \mathbb C$ as
$$
P_x(\zeta):=\frac{1}{|\psi_x'(0)|_e}\langle \psi_x(\zeta)-\psi_x(0),\psi'_x(0)\rangle,
$$
which gives the coordinate of the projection of $\psi_x(\zeta)$ to the tangent line at $x$.

Note that $P'_x(0)=| \psi_x'(0)|_e$. Given a sequence $x_n\in X$, by normality of the family $\Psi$, up to taking a subsequence if necessary, we may assume that $\psi_{x_n}\to\psi\in \widehat\Psi$ locally uniformly and $\frac{\psi'_{x_n}(0)}{|\psi'_{x_n}(0)|_e}\to v$.  It follows that $P_{\psi_{x_n}}$ converges locally uniformly to the function
\begin{equation}
\label{normallimit}
P(\zeta)=\langle\psi(\zeta)-\psi(0),v\rangle.
\end{equation}
Therefore $\mathcal P=\{P_x\,:\,x\in X\}$ is a normal family. We will write  $\widehat{\mathcal P}$ for the family of all possible normal limit of sequences in $\mathcal P$.


\begin{lemma}
\label{nonconst}
Let $f$ be quasi-expanding on $X$. Then every $P\in\widehat{\mathcal P}$ is not constant. 
\end{lemma}
\begin{proof}
The map $f^{-1}$ can be written, up to a conjugation, as a finite composition of maps of the form
$
g_i(x,y)=(y,p_i(y)-a_ix),
$
where $|a_i|\neq 0$ and $p_i$ a polynomial of degree at least $2$. 

If $R$ is sufficiently large, then $g_i(V^-_R)\subset V^-_R$, for all $g_i$, where $V^-_R$ is as in \cite[Lemma 2.5]{BS1}. Therefore if $g_1(z)\in V^-_R$, we must have $z\in U^-$. It is not hard to show that every complex plane contains a point $z$ so that $g_1(z)\in V^-_R$. Therefore $K^-$ does not contain any complex plane. 


Every $P\in\widehat{\mathcal P}$ has the form \eqref{normallimit}. If $P$ were constant than the set $\psi(\mathbb C)\subset J^-$ would be a complex plane contained in $K^-$, giving a contradiction. 
\end{proof}

Given a closed (and therefore compact) set $T\subset X$ and $x\in T$ we define
\begin{align*}
\Psi_T&:=\{\psi_x\,:\,x\in T\},\\
\widehat\Psi_T&:=\{\psi\in \widehat\Psi\,:\,\psi\text{ is a normal limit of a sequence in } T\},\\
\widehat\Psi_{x,T}&:=\{\psi\in\widehat\Psi_T\,:\,\psi(0)=x\}.
\end{align*}
Furthermore we write
\begin{equation}
\label{eq:tau}
\tau_T(x):=\max_{\psi\in\widehat\Psi_{x,T}}\textrm{Ord}(\psi),\qquad \nu_T:=\max_{x\in T}\tau_T(x).
\end{equation}

\begin{lemma}
\label{injectivnessm1}
Let $f$ be quasi-expanding on $X$. Let $T\subset X$ be a closed set so that $\nu_T=1$. Then there exists $\rho>0$ so that $P_x$ is injective on $\{|\zeta|<\rho\}$ for every $x\in X$.
\end{lemma}
\begin{proof}
The family $\mathcal P_T=\{P_x\,:\,x\in X\}$ is normal. Since $P'_x(0)=|\psi'_x(0)|_e$, and since $\nu_T=1$, the value of $P'_{x}(0)$ is uniformly bounded away from $0$. Using normality, we may find $\rho,\varepsilon>0$ so that $P'_x(\{|\zeta|<\rho\})\subset \mathbb \{\mathrm{Re}\,z>\varepsilon\}$ for every $x\in T$. Given $\zeta_1,\zeta_2\in \{|\zeta|<\rho\}$ it follows that
\begin{align*}
|P_x(\zeta_1)-P_x(\zeta_2)|&=|\zeta_1-\zeta_2|\left|\int_0^1 P'_x(t(\zeta_1-\zeta_2)+\zeta_2) \,dt\right|\\
&\ge\varepsilon|\zeta_1-\zeta_2|,
\end{align*}
proving that $P_x$ is injective.
\end{proof}

We define the \textit{unstable disk of size $\rho>0$ through $x\in X$} as 
\begin{align*}
D^u_{x,\rho}:&=\{z\in W^u_x\,|\,\mathrm{dist}^u_x(z,x)<\rho\}\\
&=\psi^u_x(\{|\zeta|<\rho\}).
\end{align*}
Both sets $W^u_{x,r}$ and $D^u_{x,\rho}$ define a neighborhood of $x\in X$, which is relatively compact in $W^u_x$. This two notions are equivalent by Propositions \ref{prop:radii} and \ref{prop:radii2}.

Note that the map $P_x$ is injective on the disk $\{|\zeta|<\rho\}$ if and only if the set $D^u_{x,\rho}$ is a graph over $E^u_{x}$, the complex plane spanned by $\psi'_x(0)$. 

\begin{figure}[h!]
\centering 
\resizebox{0.8\textwidth}{!}{ 
\begin{tikzpicture}[x=1.0cm]
\clip(-6,-2.5) rectangle (6,2.5);
\draw [->] (0,0) -- (4,0);
\draw [domain=-6:6] plot(\x,{0});
\draw[dash pattern=on 3pt off 3pt, smooth,samples=100](0,0) circle (2);
\draw (0,0)-- (0,-2);
\draw[line width=2pt, smooth,samples=100,domain=-3.0:3.0] plot(\x,{0.01*(\x)^3});
\draw[smooth,samples=100,domain=-5.0:-3.0] plot(\x,{0.01*(\x)^3});
\draw[smooth,samples=100,domain=3.0:5.0] plot(\x,{0.01*(\x)^3});
\node[above] at (0,0) {$p$};
\node[below] at (-2.5,-0.37) {$D^u_{x,\varepsilon}$};
\node[below] at (-4.2,-1.26) {$D^u_{x,\rho}$};
\node[below] at (4,-.05) {$E^u_x$};
\node[right] at (0,-1.5) {$\delta$};
\end{tikzpicture}
}
 \resizebox{0.8\textwidth}{!}{ 
\begin{tikzpicture}[x=1cm]
\clip(-6,-2.5) rectangle (6,2.5);
\draw [->] (0,0) -- (4,0);
\draw [domain=-6:6] plot(\x,{0});
\draw[smooth,samples=100,domain=0.6:1.0] plot(6*\x-3*\x^2,0.5*\x^2);
\draw[smooth,samples=100,domain=0.0:1.0] plot(-3+12*\x-6*\x^2,-\x^2+1.5);
\draw[smooth,samples=100,domain=-5.0:-3.0] plot(\x,{0.05*((\x)+3)^2+1.5});
\draw[dash pattern=on 5pt off 2pt, smooth,samples=100,domain=0:1.0] plot(2*\x^3+3,0.5*\x+0.5);
\draw[dash pattern=on 5pt off 2pt, smooth,samples=100,domain=-1.0:0] plot(-2*\x^3+3,0.5*\x+0.5);
\draw[dash pattern=on 3pt off 3pt, smooth,samples=100](0,0) circle (2);
\draw (0,0)-- (0,-2);
\draw[line width=2pt, smooth,samples=100,domain=0.0:0.6] plot(6*\x-3*\x^2,0.5*\x^2);
\draw[smooth,samples=100,domain=-5.0:-3.0] plot(\x,{0.05*(\x)^2});
\draw[line width=2pt, smooth,samples=100,domain=-3.0:0.0] plot(\x,{0.05*(\x)^2});
\node[above] at (0,0) {$p$};
\node[right] at (0,-1.5) {$\delta$};
\node[above] at (-2.5,0.4) {$D^u_{x,\varepsilon}$};
\node[above] at (-4.2,1.5) {$D^u_{x,\rho}$};
\node[below] at (4,-.05) {$E^u_x$};
\end{tikzpicture}
}
\caption{In the first picture $P_x$ is injective on $\{|\zeta|<\rho\}$, in the second it is not.}
\end{figure}
\begin{lemma}
\label{Kobelemma}
Let $f$ be quasi-expanding on $X$. Let $T\subset X$ be a closed set so that $P_{x}$ is injective on $\{|\zeta|<\rho\}$ for every $x\in T$. Then for every $0<\varepsilon<\rho$ we can find $\delta>0$ so that
\begin{equation}
\label{inclusionstableman}
D^u_{x,\rho}\cap B(x,\delta)\subset D^u_{x,\varepsilon},\qquad\forall x\in T.
\end{equation}
\end{lemma}

\begin{proof}
For every $x\in T$, the map $P_{x}$ has no critical points in $\{|\zeta|<\rho\}$. Suppose that there exists a sequence $x_n\in T$ such that $P_{x_n}'(0)\to 0$. By taking a subsequence of $x_n$ if necessary, we may assume that $\psi_{x_n}\to\psi$ and that $P_{x_n}\to P$ locally uniformly. Since $P'(0)=0$, for every $n$ sufficiently large the map $P_{x_n}$ has a critical point inside $\{|\zeta|<\varepsilon'\}$, contradicting the assumptions of the lemma. Therefore there exists $M>0$ so that $| P_x'(0)|>M$ for every $x\in T$.

For $x\in T$ the function $P_x$ is univalent in $\{|\zeta|<\varepsilon\}$. If we write $\delta=\frac{M}{4}\varepsilon$, by K\"obe Quarter Theorem,  we conclude that
$$
\{|\zeta|<\delta\}\subset P_x(\{|\zeta|<\varepsilon\}),\qquad\forall x\in T.
$$

Let $x\in T$ and $y\in D^u_{x,\rho}\cap B(x,\delta)$. Notice that $|P_x(\zeta)|<{\rm dist}(\psi(\zeta),x)$, therefore if $\zeta_y=\psi_x^{-1}(y)\in\{|\zeta|<\rho\}$, then $|P_x(\zeta_y)|< \delta$. Since $P_x$ is injective on $\{|\zeta|<\rho\}$, we must have $|\zeta_y|<\varepsilon$, and therefore $y\in D^u_{x,\varepsilon}$.
\end{proof}

\section{Characterization of hyperbolicity}\label{sec:char} 
This section is devoted to a proof of Theorem \ref{thm:2}. If a map $f$ is uniformly hyperbolic on $J$, then there are no tangencies in $J=J^*$. We will therefore assume that $f$ is quasi-hyperbolic on $J^*$ but not uniformly hyperbolic on $J$ and prove that there must be tangencies in $J^*$. In this section we will use the subscripts $s/u$ in order to distinguish between quasi-contraction and quasi-expansion. The proof of the following fact is essentially contained in \cite{BS8}.

\begin{proposition}
\label{hyperbolic}
Let $f$ be quasi-hyperbolic on $X$. Given a compact and $f$-invariant set $T\subset X$, let $\tau^{s/t}_T$ be as in \eqref{eq:tau}. Then $f$ is uniformly hyperbolic on $T$ if and only if $\nu^s_T=\nu^u_T=1$. 
\end{proposition}

By \cite{D1} the map $f$ is not uniformly-hyperbolic on $J^*$, otherwise we would have $J=J^*$. By the proposition above we either have $\nu^s:=\nu^s_{J^*}> 1$ or $\nu^u:=\nu^u_{J^*}> 1$. By replacing $f$ with $f^{-1}$ if necessary, we may always assume that $\nu^u>1$.

We define the compact and invariant set
$$
\mathcal J:=\{x\in J^*\,:\,\tau^u(x)=\nu^u\}=\{x\in J^*\,:\,\tau^u(x)\text{ is maximal}\}.
$$
We claim that $\nu^u_{\mathcal J}=1$. If instead $\nu^u_{\mathcal J}>1$, we could find a sequence $x_k\in \mathcal J$ which converges to a point $x\in \mathcal J$ with $\nu^u_{\mathcal J}$-folding along the unstable direction. For every $x_k$ we can then find a sequence $x_{k,n}$ which converges to $x_k$ with $\nu^u$-folding along the unstable direction. This implies that $$\tau^u(x)\ge \nu^s_{\mathcal J}\nu^u>\nu^u$$ contradicting the maximality of $\nu^u$.

Depending on the value of $\nu^s_{\mathcal J}$, we distinguish between two different cases
\begin{itemize}
\item[]\textbf{Hyperbolic case ($\nu^s_{\mathcal J}=1$).} In this case by Proposition \ref{hyperbolic} the set $\mathcal J$ is uniformly hyperbolic.
\item[]\textbf{Non-hyperbolic case ($\nu^s_{\mathcal J}>1$).} In this case the set $\mathcal J':=\{x\in \mathcal J\,|\,\tau^s_{\mathcal J}(x)=\nu^s_{\mathcal J}\}$ is uniformly hyperbolic although $\mathcal J$ itself is not. 
\end{itemize}

Our proof of Theorem \ref{thm:2} is easier in the non-hyperbolic case, so we will handle it first.

\subsection{Proof of Theorem \ref{thm:2} in the non-hyperbolic case}\mbox{}
\par\smallskip
\noindent
Given $y\in \mathcal J'$, by uniform hyperbolicity of $\mathcal J'$ the manifolds $W^u_y$ and $W^s_y$ intersect transversally at $y$. We choose new coordinates in a neighborhood $\mathcal N\cong \mathbb D^2$ of $y$, so that $W_{y,\mathcal N}^s=\mathbb D\times\{0\}$ and $W_{y,\mathcal N}^u=\{0\}\times \mathbb D$. Here $W^{s/u}_{x,\mathcal N}$ denotes the connected component of $W^{s/u}_x\cap\mathcal N$ containing $x\in X\cap\mathcal N$. 

Let $x_k\in \mathcal J$ be a sequence converging to $y\in \mathcal J'$ with $\nu^s_{\mathcal J}$-folding along the stable direction. Since $\nu^u_{\mathcal J}=1$, the sequence $x_k$ converges to $y$ with 1-folding along the unstable direction. Write $\pi_1,\pi_2:\mathcal N\rightarrow \mathbb D$ for the projections on the first and second coordinates of $\mathcal N$.

If $\mathcal N$ is sufficiently small and $k$ is sufficiently large, then the projection $\pi_2:W^u_{x_k,\mathcal N}\rightarrow \mathbb D$ is a degree $1$ covering map, while the projection $\pi_1:W^s_{x_k,\mathcal N}\rightarrow \mathbb D$ is a degree $\nu^s_{\mathcal J}$ branched covering map. Indeed, should there be no such $\mathcal N$ then one could find a subsequence of $x_k$ with strictly higher order of folding in one of the directions, giving a contradiction. We conclude that $ W^s_{x_k,\mathcal N}$ and $W^u_{x_k,\mathcal N}$ intersect in exactly $\nu^s_{\mathcal J}$ points counted with multiplicity.

\begin{figure}[h]
\centering 
\resizebox{0.8\textwidth}{!}{ 
\begin{tikzpicture}[x=1.0cm]
\clip(-5.1,-3.1) rectangle (5.1,3.1);
\draw(-5,-3)--(5,-3)--(5,3)--(-5,3)--(-5,-3);
\draw[dash pattern=on 1pt off 5pt](0,-3)--(0,3);
\draw[dash pattern=on 1pt off 5pt](-5,0)--(5,0);
\draw[smooth,samples=100,domain=-3:3] plot({0.02*(\x)^2+0.5},\x);
\draw[smooth,samples=100,domain=0.1:0.9] plot({30*(\x)^2-30*(\x)+7.7},\x);
\draw[dash pattern=on 5pt off 2pt, smooth,samples=100,domain=0.08367:0.9163] plot({-30*\x^2+30*\x-7.3},\x);
\node[below left] at (0,0) {$y$};
\fill  (0.5032,0.3995) circle (1.5pt);
\draw(0.72,0.17) node {\small $x_k$};
\fill  (0.5072,0.6012) circle (1.5pt);
\draw (1,0.72) node {\small $x_{k,2}$};
\node[right] at (.8,2.5){$W^u_{x_k,\mathcal N}$};
\node[below] at (4,1.5){$W^s_{x_k,\mathcal N}$};
\node[below right] at (-5,3){$\mathcal N$};
\end{tikzpicture}
}
\caption{The intersection between $W^s_{x_k,\mathcal N}$ and $W^u_{x_k,\mathcal N}$}
\end{figure}

Since $\nu^u_\mathcal J=1$, by Lemma \ref{injectivnessm1} there exists $\rho_0>0$ so that $P^u_x$ is injective on $\{|\zeta|<\rho_0\}$ for all $x\in \mathcal J$. By choosing $0<\rho<\rho_0$ small and by taking a subsequence of $x_k$ if necessary, we may assume that $D^{s/u}_{x_k,\rho}\subset W^{s/u}_{x_k,\mathcal N}$. Recall that $D^{s/u}_{x_k,\rho}=\psi^{s/u}_{x_k}(|\zeta|<\rho)$.

Given $\rho>0$, we may take $k$ sufficiently large the set $D^s_{x_k,\rho}\cap D^u_{x_k,\rho}$ contains (at least) $\nu^s_{\mathcal J}$ points counted with multiplicity. If $0<\rho<\rho_0$ is sufficiently small we conclude that
$$
D^s_{x_k,\rho}\cap D^u_{x_k,\rho}=W^s_{x_k,\mathcal N}\cap W^u_{x_k,\mathcal N}.
$$
This also shows that points in $W^s_{x_k,\mathcal N}\cap W^u_{x_k,\mathcal N}$ get uniformly close to $x_k$ with respect to the intrinsic distances $\mathrm{dist}^{s/u}_{x_k}$. We obtain the following
\begin{lemma}
\label{lemma:bula}
Let $y\in\mathcal J$ and $\rho<\rho_0$  sufficiently small. Then there exists a sequence $x_k\in\mathcal J$ converging to $y$ with $\nu^s_{\mathcal J}$-folding along the stable direction, so that $D^s_{x_k,\rho}\cap D^u_{x_k,\rho}$ contains $\nu^s_{\mathcal J}$ points counted with multiplicity, and so that
\begin{equation}
\label{convergenceintrinsic}
\lim_{k\to\infty}\mathrm{dist}^{s/u}_{x_k}\left(x_k,D^s_{x_k,\rho}\cap D^u_{x_k,\rho}\right)=0.
\end{equation}
\end{lemma}

The following proposition implies Theorem \ref{thm:2} in the non-hyperbolic case.
\begin{proposition}
\label{tangenciesnonhyp}
For some $k\ge 0$, the point $x_k$ is a tangency of order $\nu^s_{\mathcal J}$.
\end{proposition}
\begin{proof}
It suffices to find a $k$ so that the set $D^s_{x_k,\rho}\cap D^u_{x_k,\rho}$ contains a unique point. Suppose instead that for every $k$ there exists $\widetilde{x}_k\in D^s_{x_k,\rho}\cap D^u_{x_k,\rho}$ which is not $x_k$. By Proposition \ref{expansion2} the value of $\textrm{dist}^u_{f^n(x_k)}(f^n(x_k),f^n(\widetilde{x}_k))$ increases exponentially in $n$. Therefore for each $k$ there exists a unique $n_k$ so that
\begin{align*}
f^{n}(\widetilde{x}_k)&\in D^u_{f^{n}(x_k)),\rho}\qquad \forall n\le n_k\\
f^{n_k+1}(\widetilde{x}_k)&\not\in D^u_{f^{n_k+1}(x_k),\rho}.
\end{align*}
By \eqref{convergenceintrinsic} we have $\lim_{k\to\infty}{\rm dist}^u_{x_k}(x_k,\widetilde x_k)=0$. Therefore the sequence $n_k$ is divergent.

Let  $\varepsilon\in (0,\rho)$ small enough so that $f(D^u_{x,\varepsilon})\subset D^u_{f(x),\rho}$ for every $x\in J^*$, and let $\delta>0$ as in Lemma \ref{Kobelemma}.

As a consequence of Propositions \ref{prop:radii} and \ref{expansion2}, there exists $n_\delta>0$ so that
$$
f^{n}(D^s_{x,\rho})\subset W^s_{f^{n}(x),\delta}\subset B(f^{n}(x),\delta),\qquad\forall n\ge n_\delta,\,x\in J^*.
$$
Choose $k$ sufficiently large so that $n_k\ge n_\delta$. By invariance of $\mathcal J$ and \eqref{inclusionstableman} we conclude that
$$
f^{n_k}(\widetilde{x}_k)\in  D^u_{f^{n_k}(x_k),\rho}\cap B(f^{n_k}(x_k),\delta)\subset D^u_\varepsilon(f^{n_k}(x_k)).
$$
By the definition of $\varepsilon>0$ it follows that $f^{n_k+1}(\widetilde{x}_k)\in D^u_{f^{n_k+1},\rho}$.
This contradicts the definition of $n_k$, showing that tangencies exist.
\end{proof}

\subsection{Proof of Theorem \ref{thm:2} in the hyperbolic case} \mbox{}
\par\smallskip
\noindent
In the non-hyperbolic case we were able to find tangencies inside $\mathcal J$, using that $D^u_{x,\rho}$ are graphs over the unstable direction. When $\mathcal J$ is hyperbolic, we will find our tangencies inside $J^*\setminus \mathcal J$. But now by Lemma \ref{injectivnessm1} there exists $\rho_0>0$ such that $P^{s/u}_x$ is injective on $\{|\zeta|<\rho_0\}$ for every $x\in\mathcal J$.

We will use the following analogue of Lemma \ref{lemma:bula} for the hyperbolic case.
\begin{lemma}
Let $y\in \mathcal J$ and $\rho<\rho_0$ sufficiently small. Then there exists a sequence $x_k\in J^*\cap W^s_y$ which converges to $y$ with $\nu^u$-folding along the unstable direction, so that $\lim_{k\to\infty}{\rm dist}^s_y(y,x_k)=0$, so that  $D^s_{x_k,\rho}\cap D^u_{x_k,\rho}$ contains $\nu^u$ points counted with multiplicity, and so that
$$
\lim_{k\to\infty}\mathrm{dist}^{s/u}_{x_k}\left(x_k,D^s_{x_k,\rho}\cap D^u_{x_k,\rho}\right)=0.
$$
\end{lemma}
\begin{proof}
Choose a sequence $x_k\in J^*$ which converges to $y\in \mathcal J$ with $\nu^u$-folding along the unstable direction. 

Similarly to the non-hyperbolic case, we can find a neighborhood $\mathcal N$ of $y$ so that when $k$ is sufficiently large $W^s_{y,\mathcal N}\cap W^u_{x_k,\mathcal N}$ contains $\nu^u$ points counted with multiplicity. Given $k$ large, we choose $\widetilde x_k\in W^s_{y,\mathcal N}\cap W^u_{x_k,\mathcal N}$. As in the previous section, we may show that the points in the sequence $\widetilde x_k$ get close to $y$ and $x_k$ with respect to the intrinsic distances $\mathrm{dist}^s_{y}$ and $\mathrm{dist}^u_{x_k}$ as $k\to\infty$. This shows that $\lim_{k\to\infty}{\rm dist}^s_y(y,\widetilde x_k)=0$.

For every $k$ we have $\psi^u_{\widetilde x_k}(\zeta)=\psi^u_{x_k}(a_k\zeta+b_k)$. Since $\mathrm{dist}^u_{x_k}(x_k,\widetilde x_k)$ converges to zero, by Lemma \ref{intrinsicmetric} we conclude that $|a_k|\to 1$ and that $b_k\to 0$. This shows that the sequences $\psi^u_{\widetilde x_k}$ and $\psi^u_{x_k}$ have the same normal limit, and thus that $\widetilde x_k$ converges to $y$ with $\nu^u$-folding along the unstable direction. Furthermore we have 
$$
W^s_{\widetilde x_k,\mathcal N}=W^s_{y,\mathcal N}\quad\text{and}\quad W^u_{\widetilde x_k,\mathcal N}=W^u_{x_k,\mathcal N},
$$
showing that the set $W^s_{\widetilde x_k,\mathcal N}\cap W^u_{\widetilde x_k,\mathcal N}$ consists of $\nu^u$ points counted with multiplicity. Following the strategy used in the non-hyperbolic case, it is not hard to prove that the sequence $\widetilde x_k$ satisfies the conclusions of the Lemma.


\end{proof}

In the remainder of this section we will assume that $f$ has no tangencies on $J^*$, which in the end will lead to a contradiction. Let $y\in \mathcal J$ and $0<\rho<\rho_0$ small, and choose a sequence $x_k$ as in the previous lemma. The absence of tangencies implies that the set $D^s_{x_k,\rho}\cap D^u_{x_k,\rho}$ consists of $\nu^u$ distinct points, which we will denote as $\{x_k,x_{k,2},\dots,x_{k,\nu^u}\}$. By changing their order if necessary, we may assume that there exists a sequence $m_k$, so that $\lim_{k\to\infty} m_k=\infty$ and so that
\begin{align*}
f^{-n}(x_{k,i})&\in D^s_{f^{-n}(x_k),\rho}\qquad \forall n\le m_k,\,\forall i=2,\dots, \nu^u\\
f^{-m_k-1}(x_{k,2})&\not\in D^s_{f^{-m_k-1}(x_k),\rho}.
\end{align*}

For every couple of positive integers $(k,n)$ we write $\psi_{k,n}=\psi^s_{f^{-n}(x_k)}$ and $P_{k,n}=P^s_{f^{-n}(x_k)}$. Note that if for some large $k$ the map $P_{k,m_k}$ were injective, then the strategy of Proposition \ref{tangenciesnonhyp} would work, giving existence of tangencies. Since we are under the assumption that tangencies do not exist, it follows that
\begin{lemma}
If $k$ is sufficiently large, the map $P_{k,m_k}$ is not injective on $\{|\zeta|<\rho\}$.
\end{lemma} 


After taking a subsequence of $x_k$ if necessary, we may assume that all maps $P_{k,m_k}$ are not injective on $\{|\zeta|<\rho\}$. Since $\lim_{k\to\infty}\mathrm{dist}^s_y(y,x_k)=0$, by Lemma \ref{intrinsicmetric}, we have that $\mathrm{dist}^s_{x_k}=c_k\,\mathrm{dist}^s_y$, with $c_k\to1$. It follows that when $k$ is sufficiently large $D_{x_k,\rho}^s\subset D^s_{y,\rho_0}$, proving that the map $P_{k,0}$ is injective on $\{|\zeta|<\rho\}$.

Therefore we may find a positive integer $0\le n_k< m_k$ so that
\begin{align*}
&P_{k,n}\textrm{ is injective on }\{|\zeta|<\rho\}\textrm{ for all }n\le n_k,\\
&P_{k,n_k+1}\textrm{ is not injective on }\{|\zeta|<\rho\}.
\end{align*}
The sequence $x_k$ converges to $y$ with respect to the intrinsic distance $\mathrm{dist}^s_y$. Therefore for every $n\ge 0$ we have $D^s_{f^{-n}(x_k),\rho}\subset D^s_{f^{-n}(y),\rho_0}$ for every $k$ sufficiently large. Therefore the sequence $n_k$ diverges.

By taking a subsequence of $x_k$ if necessary, we may assume that for every integer $i$ the sequence $\psi_{k,n_k+i}$ converges locally uniformly to $\psi_i\in\widehat\Psi^s$ and $P_{k,n_k+i}$ converges locally uniformly to $P_i\in\widehat{\mathcal P}^s$. Furthermore we have
$$
P_{i}(\zeta)=\langle\psi_i(\zeta)-\psi_i(0),v_i\rangle,
$$
for some unit vector $v_i$.

As shown in \cite[Lemma 3.2]{GP}, for every integer $i$ there exists $\lambda_i\in \mathbb C$ so that
\begin{equation}
\label{transitionformula}
f^{-i}\circ\psi_0(\zeta)=\psi_i(\lambda_i\,\zeta).
\end{equation}
Furthermore if $\kappa>1$ is as in Proposition \ref{prop:eqcond} for the map $f^{-1}$, then we have $|\lambda_i|\ge\kappa^{i}$ when $i>0$ and $|\lambda_i|\le\kappa^{i}$ when $i<0$.
\begin{lemma}
\label{injectivepar}
The maps $\psi_i:\mathbb C\rightarrow\mathbb C^2$ are injective holomorphic immersions.
\end{lemma}
\begin{proof}
Suppose that the map $\psi_0$ has a critical point at $\zeta_0$. Let $\kappa>1$ as above and choose a positive integer $N$ so that $\kappa^{-N}\zeta_0\in\{|\zeta|<\rho\}$. By \eqref{transitionformula} we have
$$
f^{N}\psi_0(\zeta)=\psi_{-N}(\lambda_{-N}\,\zeta)
$$
with $|\lambda_{-N}|\le\kappa^{-N}$. We conclude that the map $\psi_{-N}$ has a critical point at $\zeta_{-N}=\lambda_{-N}\zeta_0\in\{|\zeta|<\rho\}$. It is clear that $P_{-N}'(\zeta_{-N})=0$. 

By Lemma \ref{nonconst} the map $P_{-N}$ is not constant,  therefore $\zeta_{-N}$ is an isolated critical point. By Rouch\'e Theorem, we deduce that when $k$ is sufficiently large $P_{k,n_k-N}$ also has a critical point close to $\zeta_{-N}$ and inside $\{|\zeta|<\rho\}$, contradicting the fact that $P_{k,n}$ is injective on this disk for $n\le n_k$.

Write $y_0=\psi_0(0)$. Then by \cite[Corollary 3.1]{GP}, we have that $\psi_0=\psi_{y_0}\circ h$ for some polynomial $h$.  Since $\psi_0$ has no critical point on all $\mathbb C$ we conclude that $h$ is an affine map proving that $\psi_0$ is an holomorphic immersion. By \eqref{transitionformula}, given an integer $i$ we have $\psi_i(\zeta)=f^{-i}\psi_0\left(\zeta\,\lambda_i^{-1}\right)$, proving that every $\psi_i$ is an injective holomorphic immersion.
\end{proof}

If $y_i=\psi_i(0) $, as a consequence of Proposition \ref{prop:bigfam}, we must have that $\psi_i$ and $\psi_{y_i}$ are equal up to a rotation. It is also not hard to show that 
$$
P_i(\zeta)=\frac{1}{|\psi'_i(0)|_e}\langle\psi_i(\zeta)-\psi_i(0),\psi_i'(0)\rangle
$$

Since $n_k< m_k$, given $\kappa>1$ as in Proposition \ref{expansion2} for the map $f^{-1}$, we have that 
$$
f^{-n_k}\{x_k,x_{k,2},\dots,x_{k,\nu^u}\}\subset D^u_{f^{-n_k}(x_k),\kappa^{-n_k}\,\rho}\cap D^s_{f^{-n_k}(x_k),\rho}.
$$

Recall that $\psi_0=\lim_{k\to\infty}\psi^s_{f^{-n_k}(x_k)}$. Lemma \ref{injectivepar} implies that the sequence $f^{-n_k}(x_k)$ converges to $y_0=\psi_0(0)$ with 1-folding along the stable direction. Every intersection point $f^{-n_k}(x_{k,i})$ gets uniformly close to $f^{-n_k}(x_k)$ with respect to the intrinsic metric ${\rm dist}^u_{f^{-n_k}(x_k)}$. Therefore by continuity of transverse intersection either $y_0$ is a tangency point, which is not possible by the no tangencies assumption, or the sequence $f^{-n_k}(x_k)$ converges to $y_0$ with $\nu^u$-folding along the unstable direction, which implies that $y_0\in \mathcal J$. Since the latter must be the case we deduce that $f^{-1}(y_0)=\psi_1(0)\in \mathcal J$.

Since $f^{-1}(y_0)\in \mathcal J$  we conclude that $P_{1}$ is injective on $\{|\zeta|<\rho_0\}$ and that $P_{k,n_k+1}$ converges locally uniformly to this function. It follows that for $k$ sufficiently large, the function $P_{k,n_k+1}$ is injective on $\{|\zeta|<\rho\}$, contradicting the definition of $n_k$, and therefore the assumption that there are no tangencies in $J^*$. This concludes the proof of the Theorem.

\section{Consequences of Theorem \ref{thm:2}}

In this section we will assume that $f$ is quasi-hyperbolic on $J^*$. Given $m^s,m^u\ge 1$ we define
$$J^*_{m^s,m^u} := \{x\in J^*:  \tau^s(x)=m^s, \tau^u(x)=m^u\}.$$


\begin{proposition}
\label{prop:tang}
If a point $x\in J^*_{m^s,m^u}$ and $W^u_x$ and $W^s_x$ are tangent to order $k$ then the forward limit set of $x$ is contained in $J^*_{p,q}$ with $p\ge m^s$ and $q\ge (k+1)m^u$ the backward limit set of $x$ is contained in $J^*_{p',q'}$ with $p'\ge (k+1)m^s$ and $q'\ge m^u$.
\end{proposition}

\begin{proof}
We prove the first statement. The second statement follows by considering $f^{-1}$. Consider the case when $m^s=m^u=1$. Since $W^u_x$ and $W^s_x$ are tangent to order $k$ we can assume that for some $c\ne0$, $\psi^s_x(\zeta)=\psi^u_x(c\zeta)+\cdots$ where the dots represent terms or order greater than $k$. Write $f^n(\psi^s_x(\zeta))=\psi^s_{f^n(x)}(\lambda_n\zeta)$ and $f^n(\psi^u_x(\zeta))=\psi^u_{f^n(x)}(\mu_n\zeta)$. For some $\kappa>1$ we have $\lambda_n\ge\kappa^n$ and $\mu_n\le\kappa^{-n}$. 

Now $f^n(\psi^s_x(\zeta))=f^n(\psi^u_x(c\zeta))$ up to terms of order greater than $k$. Write $\psi^s_{f^n(x)}(\zeta)=f^n(x)+\vec a_{1,n}\zeta+\vec a_{2,n}\zeta^2+\cdots$ and $\psi^u_{f^n(x)}(\zeta)=f^n(x)+\vec b_{1,n}\zeta+\vec b_{2,n}\zeta^2+\cdots$. Then we have $\psi^s_{f^n(x)}(\lambda_n\zeta)$ is equal to $\psi^u_{f^n(x)}(c\mu_n\zeta)$ up to order $k$ so by comparing coefficients we get $\vec a_{j,n}\lambda_n^j=\vec b_{j,n}c^j\mu_n^j$ for $j=1,\dots, k$. By normality, $|b_{j,n}|$ is bounded by a constant depending on $j$ but not $n$ so 
$|\vec a_{j,n}|\le C_j \kappa^{-2n} |c^j|$. We conclude that for $j=1,\dots ,k$, $|\vec a_{j,n}|\to 0$ as $n\to\infty$. In particular if $x'$ is in the forward limit set of $x$ then taking $\psi\in\widehat\Psi_{x'}^u$ to be a convergent subsequence of $\psi^u_{f^n(x)}$ the coefficients of $\psi(\zeta)$ vanish for $j=1,\dots, k$.

If there is a map $\psi_0\in\widehat\Psi_x^u$ which vanishes to order $\ell$ then $\psi$ can be written as the composition of the map $\psi^u_x$ and a map $\alpha$ of order $\ell$. Arguing as above we can find a sequence $\psi_n\in\widehat\Psi_{f^n(x)}^u$ whose limit vanish to order at least $(k+1)\ell$.
\end{proof}

\begin{corollary}
There is a bound on the order of tangency between stable and unstable manifolds.
\end{corollary}
\begin{proof} 
In \cite{BS8} it is shown that the set $J^*_{m^s,m^u}$ is empty for $m^s$ and $m^u$ large. 
\end{proof}

\begin{corollary}If $q$ is a point of tangency between stable and unstable manifolds, then $q$ is wandering, i.e., it is not contained in either its forward or backward limit set.
\end{corollary}

\begin{proof}
We have $q\in J^*_{m^s,m^u}$ for some $m^s$ and $m^u$.  By Proposition \ref{prop:tang}, every point of the forward or backward limit set of a point of tangency is contained in $J^*_{n^s,n^u}$, with either $n^s>m^s$ or $n^u>m^u$, or both. 
\end{proof}

\begin{theorem}
If $f$ is quasi-hyperbolic on $J^*$, then the condition that $J^+$ is laminated in a neighborhood of $J^*$ is equivalent to the condition that $J^-$ is laminated a neighborhood of $J^*$; and either condition is equivalent to hyperbolicity on $J$. 
\end{theorem}

\begin{proof} 
If $f$ is quasi-hyperbolic on $J^*$ but not hyperbolic on $J$, then by Theorem \ref{thm:2} there is a tangency of order $k\ge 2$. Therefore by Proposition \ref{prop:tang} both $J^*_{j_s,k}$ and $J^*_{k,j_u}$ are non-empty for some $j_s,j_u\ge1$.  It follows that the disks of $W^s_{x,r}$ (as well as the disks of $W^u_{x,r}$) exhibit local folding at some points of $J^*_{\rm max}$.   If $p$ is a saddle point, then the stable manifold $W^s_p$ is dense in $J^+$, and each disk $W^s_{x,r}$ is either disjoint from or contained in $W^s_p$.  If $J^+$ is laminated, then the disks $W^s_{x,r}$ are subsets of global leaves of the lamination.  In conclusion, if $f$ is not uniformly hyperbolic, then the $W^s_{x,r}$ exhibit local folding at some point $x_0$, and so $J^+$ is not laminated in a neighborhood of $x_0$.
\end{proof}

Let $L$ be a Riemann Surface injectively immersed in $\mathbb C^2$. Given $x\in L$ we write $L_{x,\varepsilon}$ for the connected component of $L\cap B(x,\varepsilon)$ containing $x$. Here the connected component is taken with respect to the intrinsic topology of $L$. Given another Riemann surface $S\subset \mathbb C^2$  and $x\in L\cap S$, we say that $L$ and $S$ \emph{are locally equal at } $x$ if $L_{x,\varepsilon}=S_{x,\varepsilon}$ for some $\varepsilon>0$. If $L$ and $S$ are not locally equal at $x$, then there exists $\varepsilon>0$ so that $L_{x,\varepsilon}\cap S_{x,\varepsilon}=\{x\}$.

Finally, we recapture an earlier result from \S2 of \cite{BedD}. The following theorem provides a characterization of hyperbolicity that does not assume a priori that $f$ is quasi-hyperbolic, and it removes the transversality hypothesis from \cite[Theorem 8.3]{BS8}:

\begin{theorem}
\label{thm:lam}
A H\'enon map $f$ is uniformly hyperbolic on $J$ if and only if $J^+$ and $J^-$ are Riemann surface laminations in a neighborhood of $J^*$.
\end{theorem}
\begin{proof}
By the previous theorem, it suffices to show that if $J^\pm$ are laminated in a neighborhood of $J^*$, then $f$ is quasi-hyperbolic on $J^*$. We will write $L^{s}_x\in \mathcal L^+$ and $L^u_x\in\mathcal L^-$ for the \emph{stable and unstable} leaves of the lamination through a point $x\in J^*$.

Assume that $S\subset K^+$ is a Riemann surface containing $x\in J^*$. We claim that $L^s_x$ and $S$ are locally equal at $x$. By Slodkowski \cite{S} there exists an open neighborhood $\mathcal U\ni x$ so that $\mathcal L^+|_\mathcal{U}$ extends to a lamination $\mathcal L^*$ of $\mathcal U$. Every leaf of $\mathcal L^*$ is either contained in $K^+$ or in $U^+$. Assume that $L^s_x$ and $S$ are not locally equal at $x$. Then for $y\in \mathcal U$ close to $x$, the leaf of $\mathcal L^*$ through $y$ intersects $S$, and therefore is contained in $K^+$. We conclude that a ball centered at $x$ is contained in $K^+$, which is not possible since $x\in J^*$.

If $p$ is a saddle then $L^s_p\subset W^s_p$ and $L^u_p\subset W^u_p$. In particular the stable and unstable manifolds satisfy the proper bounded area condition of \cite{BS8}.
By \cite[Corollary 3.5]{BS8} the map $f$ is quasi-expanding provided that for every $\delta>0$ there exists $\eta>0$ so that 
$$
\sup_{W^u_{p,\delta}} G^+>\eta\qquad\text{for every saddle }p\in\mathcal S.$$A similar statement holds for quasi-contraction. If the condition above is not satisfied then, since the leaves of the lamination change continuously, the leaf $L^u_x\in \mathcal L^-$ is locally contained in $K^+$ for some $x\in J^*$. In particular the stable and unstable leaves $L^{s/u}_x$ are locally equal at $x$. The theorem follows once we show that this cannot occur.

The laminations $\mathcal L^\pm$ are locally invariant, meaning that for $x\in J^*$ the images of $L^{s/u}_x$ under the map $f$ are locally equal to $L^{s/u}_{f(x)}$ at $f(x)$. This holds for saddle points, and thus it holds everywhere in $J^*$ by continuity of the lamination. It follows that $\mathcal L^\pm$ extend to laminations of $J^\pm$. Indeed given a point $y\in J^+$, by the Ergodic Closing Lemma of \cite{D2} every cluster value of the sequence of Cesar\`o averages $\nu_n=n^{-1}\sum_n\delta_{f^n(y)}$ is supported on $J^*$. In particular there exists a subsequence $n_k$ so that $f^{n_k}(y)\to J^*$. By pulling back the lamination $\mathcal L^+$ along the orbit of $y$, we obtain a local lamination near $y$. The  invariance of $\mathcal L^+$ guarantees that the lamination near $y$ is locally unique. We can glue together these local laminations to obtain a global lamination of $J^+$. We will keep denoting the global laminations of $J^\pm$ as $\mathcal L^\pm$.

Suppose that $L^s_x$ and $L^u_x$ are locally equal at $x\in J$. If $L^s_x\not\subset L^u_x$, we write $V$ for the connected component of $\mathrm{int}(L^s_x\cap L^u_x)$ containing $x$, and we choose $y\in\partial V$. It is important to remark that the topology we consider here is the intrinsic topology of the Riemann surface $L^s_x$. Since $y\in J$, stable and unstable leaves are defined at $y$, and we have $L^s_y=L^s_x$. Choose a neighborhood $\mathcal U\ni y$ such that $\mathcal L^-|_\mathcal{U}$ is homeomorphic to the trivial lamination of $S\times\mathbb D$.  Let $V'$ be the connected component of $V\cap \mathcal U$ so that $y\in \partial V'$. The set $V'\subset J^-$ is a Riemann surface, thus $V'$ is contained in a single leaf of the lamination $\mathcal L^-|_\mathcal{U}$. In particular $V'\subset L^u_y$, and therefore $L^u_y=L^u_x$. 

It follows that the local intersection between $L^s_x$ and $L^u_x$ at $y$ always contains more than one point, and therefore that $L^s_x$ and $L^u_x$ are locally equal at $y$. This contradicts the fact that $y\in \partial V$. We conclude that $L^s_x\subset L^u_x$, and similarly that $L^u_x\subset L^s_x$, which implies that $L^s_x=L^u_x$.  Given $y_0\in \overline {L^s_x}$ we have $y_0\in J^+$. By the maximum principle, the point $y_0$ is not a local maximum of the map $y\mapsto{\rm dist}(y,0)$ on the manifold $L^s_{y_0}$. Since the leaves of $\mathcal L^+$ move continuously, there exists a point $y'\in L^s_x$ with ${\rm dist}(y',0)>{\rm dist}(y_0,0)$. This shows that $L^s_x$ is not bounded, contradicting the fact that $L^s_x=L^u_x\subset J$, and thus showing that $L^s_x$ and $L^u_x$ cannot be locally equal at any point of $J$, concluding the proof of the theorem.
\end{proof} 


\printbibliography
\end{document}